\documentstyle[amssymb,12pt]{article}
\begin{document}
\newtheorem{guess}{Proposition }[section]
\newtheorem{theorem}[guess]{Theorem}
\newtheorem{lemma}[guess]{Lemma}
\newtheorem{corollary}[guess]{Corollary}
\newtheorem{remark}[guess]{Remark}
\def\square#1{\vbox{\hrule
\hbox{\vrule\hbox to #1 pt{\hfill}\vbox{\vskip #1 pt}\vrule}\hrule}}
 
 \centerline{\em \LARGE \bf Two applications of instanton numbers}
\vspace{8mm}

\centerline{\Large Elizabeth Gasparim \footnote{ Supported by the 
Isaac Newton Institute and  NMSU summer research award.}}
\vspace{7mm}
 
\begin{abstract}{ The two applications are: 1. sometimes
instanton numbers stratify moduli of bundles better than Chern numbers.
2. sometimes  instanton numbers 
distinguish  singularities better than the classical 
numerical invariants.  }
\end{abstract}

\section{Introduction}In a little more detail,
the two applications are: $1'.$ instanton numbers give the coarsest 
stratification of moduli of bundles on blow-ups for which 
the strata are separated.
$2'.$ some analytically inequivalent plane curve singularities have 
same $\delta_P,$ Milnor number  and Tjurina number, but distinct 
instanton numbers. The instanton numbers we use are 
local analytic invariants for instantons on a blow-up.

Let $\widetilde{{\Bbb C}^2}$ denote the blow-up of 
${\Bbb C}^2$ at the origin.
Rank 2 instantons on 
  $\widetilde{{\Bbb C}^2}$
 are built from simple algebraic data, namely, a triple 
$(j,p,t_{\infty}),$ made by an integer $j,$ 
a polynomial $p,$ 
and a framing at infinity, that is, a holomorphic map
$t_{\infty}\colon {\Bbb C}^2-\{0\} \rightarrow
{\mbox GL}(2,{\Bbb C}).$
These instantons have two holomorphic invariants:
the height and  the width, whose 
sum gives the topological charge.
Here I give two applications of these instanton numbers.
First I  use these
numbers to stratify moduli of instantons on the blown-up
plane and second I use  
this pair 
as analytic 
invariants for  plane curve singularities.
 I show that
the pair  (height, width) gives 
instanton invariants that are strictly finer than the 
topological charge of the instanton.
In fact, the stratification
of moduli of instantons by this pair of invariants is strictly finer 
that the stratification  by topological charge.
As applications to singularities, I show that these numbers 
 distinguish nodes/tacnodes from cusps/higher order cusps.
I also give an example of analytically inequivalent
curve singularities that are not distinguished by  the 
classical invariants (Milnor number, Tjurina number and the $\delta_P$ 
invariant which calculates the change in arithmetic genus)
but have distinct instanton numbers.

The charge of a  rank 2 instanton on  
$\Delta =(j,p,t_{\infty})$ on $\widetilde{{\Bbb C}^2}$ 
ranges between $j$ and $j^2$ 
depending on $p.$ 
However, unlike instantons on $S^4,$
whose charge is given locally by a unique invariant, called  
the
multiplicity, these instantons  have two
local holomorphic invariants: the height and the width. 
These invariants do not depend on the framing, and 
neither does the topological charge.
We can therefore calculate height, width and charge 
directly from the algebraic data $(j,p).$
To use 
instanton numbers as invariants of curve singularities the trick is as follows.
Given a plane curve  $p(x,y)=0$ with singularity at the origin,
chose an  
integer $j,$  and  construct an instanton 
with data 
$(j,p).$ Then  use its 
numerical invariants  as 
analytic invariants of the curve. 

 Instantons on $\widetilde{{\Bbb C}^2}$ their moduli 
and their topological and holomorphic  invariants are described 
in section 2 and used to stratify moduli of instantons.
In section 3, 
these invariants are used to distinguish  plane curve singularities.

\vspace{5mm}
\noindent {\bf Acknowledgments}: I am grateful to  Professor
Nigel Hitchin for encouraging me to work on instanton numbers and to
Professors Robin Hartshorne, John Tate, and Abramo Hefez for enlightening 
discussions about curve singularities.
This paper was written when I participated at the 
High Dimensional Geometry program at the Isaac Newton 	Institute.
I  thank the organizers Alessio Corti, Mark Gross and 
Miles Reid
for inviting me to the program.
 During my visit, the Department of Mathematics 
of Cambridge University provided  ideal working conditions that
inspired me
to finally  write  this paper that I had long postponed.

\section{Instantons on  $\widetilde{{\Bbb C}^2}$} 

We  show that every rank 2
 instanton on $\widetilde{{\Bbb C}^2}$ is determined
by a triple $(j,p,t_{\infty}),$
where $j$ is an integer called the {\em splitting type} 
of the instanton, $p$ a polynomial and $t_{\infty}$
a trivialization at infinity.
Generically, two 
triples
$(j,p,t_{\infty})$ and 
$(j', p', t'_{\infty})$
determine the same instanton
if an only if $j'=j,$ $p' = \lambda p$
and $t'_{\infty} = A \, t_{\infty}$
where $\lambda \neq 0,$ and 
$A \in \Gamma \left({\Bbb C}^2 - \{0\},{\mbox GL}(2,{\Bbb C})\right).$ 
An instanton
$\Delta=(j,p,t_{\infty})$ 
is generic if and only if its topological charge  equals its splitting type
$j.$
Moreover,
for  every $j > 1$ there are 
nongeneric instantons $(j,p,t_{\infty}),$ with topological charge varying 
from $j+1$ up to  $j^2.$
 For each integer $j,$ we  topologize the set
${\cal M}_j$ of equivalence classes 
 of  instantons 
$(j,p,t_{\infty})$ and show that the generic set 
is a ${\mbox GL}(2,{\Bbb C})$--bundle over 
a quasi-projective smooth variety of complex dimension ${2j-3}$ 

The fact that  an instanton on  $\widetilde{{\Bbb C}^2}$
 is determined by a triple $(j,p,t_{\infty})$
follows essentially from putting together
two results: first, the proof due to King  \cite{KI}
 of the Hitchin--Kobayashi correspondence 
 over the noncompact surface  $\widetilde{{\Bbb C}^2}$
and second, the characterization of rank two holomorphic bundles on 
 $\widetilde{{\Bbb C}^2}$ given in \cite{BSPM}.
We review these two results.

Instantons on the blown-up plane are naturally identified with instantons on
 $\overline{\Bbb CP}^2$ framed at infinity; this is a simple consequence of the 
fact that  $\overline{\Bbb CP}^2$ is the conformal compactification of
$\widetilde{{\Bbb C}^2}.$ 
On his Ph.D. thesis, A. King \cite{KI}
 identifies the moduli space MI$(\widetilde{{\Bbb C}^2};r,k)$
of instantons on the blown-up plane of rank $r$ and charge $k,$ with
the moduli space
MI$(\overline{\Bbb CP}^2, \infty: r, k)$
of instantons on $\overline{\Bbb CP}^2,$
framed at $\infty,$ whose underlying vector bundle
has rank $r,$ and Chern classes $c_1 = 0$ and $c_2 = k.$

On the other hand, 
we may consider the canonical complex compactification 
of $\widetilde{\Bbb C}^2, $ which is the 
the Hirzebruch surface $\Sigma_1,$
obtained from  $\widetilde{{\Bbb C}^2}$ by adding a line ${\ell}_{\infty}$
at
infinity. 
Essentially {\em by definition}
King  identifies the moduli space
MH$(\widetilde{{\Bbb C}^2};r,k)$ of ``stable''
holomorphic bundles on $\widetilde{{\Bbb C}^2}$ 
with rank $r$ and $c_2 = k$ with the moduli space 
MH$(\Sigma_1,\ell_{\infty};r,k)$
of holomorphic bundles on $\Sigma_1$ with a trivialization 
along $\ell_{\infty}$ 
and whose underlying vector bundle has 
rank $r,$ $c_1=0$ and $c_2=k.$
King then proves the Hitchin--Kobayashi correspondence in this case, 
namely that the map $$ \mbox{MI}(\widetilde{{\Bbb C}^2};r,k) \rightarrow 
\mbox{MH}(\widetilde{{\Bbb C}^2};r,k)$$
given by taking the holomorphic part of 
an instanton connection is a bijection.
Therefore, a rank 2 instanton on 
$\widetilde{{\Bbb C}^2}$   is completely determined by 
a rank two holomorphic bundle on $\widetilde{{\Bbb C}^2}$
with vanishing first Chern class,
 together with a trivialization at infinity. 
The instanton has charge $k$ if and only if 
the corresponding holomorphic bundle 
 extends to a bundle
on $\Sigma_1$ trivial on $\ell_{\infty}$ having $c_2=k.$

We are thus led to study holomorphic rank two bundles on 
$\widetilde{{\Bbb C}^2}$ with vanishing first Chern class.
It turns out that holomorphic bundles on $\widetilde{{\Bbb C}^2}$
 are algebraic, they are extensions of line bundles
 and moreover they are trivial on 
the complement of the exceptional divisor (see \cite{CA} and \cite{JA}).
Triviality outside the exceptional divisor in this case is 
very useful and is intrisically related to the fact that 
we have algebraic bundles. It is of course not true in general that
holomorphic bundles on $\widetilde{{\Bbb C}^2} - \ell$ are trivial.

A holomorphic rank 2 bundle $E$ on $\widetilde{{\Bbb C}^2}$
with vanishing first Chern class splits over the exceptional divisor 
as ${\cal O}(j) \oplus {\cal O}(-j)$ for some positive integer $j,$
called the {\it splitting type} of the bundle, and, in this case, $E$ is
 an algebraic extension
\begin{equation}\label{extension}0 \rightarrow {\cal O}(-j) 
\rightarrow E \rightarrow {\cal O}(j) \rightarrow 0 \end{equation}
(here by abuse of  notation we write  ${\cal O}(k)$ 
 both for the line bundle ${\cal O}(k)$
 over the exceptional divisor $\ell$ as 
well as for its pull-back to  $\widetilde{{\Bbb C}^2}$).
A bundle $E$ fitting in  an exact sequence (1) 
 is determined by its extension class in 
$p \in Ext^1({\cal O}(-j), {\cal O}(j)),$ where $p$ a polynomial,
since as showed in $\cite{CA1}$ the bundle $E$  is actually  algebraic.
To this bundle on  we assign
a {\em canonical form} of  transition matrix.
We fix, once and for all, the following  charts: 
 $\widetilde{{\Bbb C}^2} = U \cup V$ where  
$U = \{(z,u)\} \simeq   {\Bbb C}^2 \simeq \{(\xi,v)\} = V$
 with $(\xi,v)=(z^{-1}, zu)$
in $U \cap V.$ 
Once these charts are fixed,  
$E$ has the canonical transition 
matrix of the form  (see
 \cite{JA} Thm.\thinspace2.1)
\begin{equation}\label{matrix}\left(\matrix {z^j & p \cr 0 &  z^{-j}} \right)\end{equation}
from $U$ to $V,$  where 
\begin{equation}\label{polynomial}p \colon = \sum_{i = 1}^{2j-2} \sum_{l = i-j+1}^{j-1}p_{il}z^lu^i \end{equation}
is a polynomial in  $z, \,z^{-1}$ and $ u.$

It follows that a rank 2 holomorphic bundle $E$ on $\widetilde{{\Bbb C}^2}$
with vanishing first Chern class is completely determined by a pair 
$(j,p)$ where $j$ is a nonnegative integer and a $p$ 
is a polynomial of the form (\ref{polynomial}). 
According to King's results, to have an instanton
we need also  a trivialization at infinity.
However, it follows  from \cite{CA} Cor.\thinspace 4.2,
that the bundle $E$ is trivial 
outside the exceptional divisor. Therefore, to any 
bundle  over $\widetilde{{\Bbb C}^2}$
represented by a pair $(t,p)$
we may assign a trivialization at infinity $t_{\infty}
\in {\mbox GL}(2, {\Bbb C}^2 -\{0\})$
thus obtaining  an instanton.
As a consequence every rank--two
 instanton  $\Delta$ on $ \widetilde{{\Bbb C}^2}$
is determined by a triple 
\begin{equation}\label{triple}\Delta:=(j,p,t_{\infty}).\end{equation}

To define the topological charge of the instanton we need to extend 
$\Delta$ to a bundle on a compact surface.
The charge is independent of the chosen compactification
(and in fact it only depends on an infinitesimal neighborhood 
of the exceptional divisor), but for simplicity
we may take the compactification of $\widetilde{{\Bbb C}^2}$
 be the Hirzebruch surface $\Sigma_1.$
This extension is obtained as follows.
Let $M$ be a complex manifold 
and $N$ a complex submanifold of $M.$
We denote by ${\cal E}_r(M,N)$ the set of equivalence 
classes of pairs $(E, \eta)$
where $E$ is a rank $r$ holomorphic bundle over $M$ 
such that $E|_N$ is trivial, and $\eta$ is a trivialization of 
$E|_N.$
Here $(E,\eta)$ is equivalent to $(E', \eta ')$
if there is a bundle equivalence $\alpha\colon E \rightarrow E'$
such that $\alpha \, \eta = \eta'.$
Recall that $\ell$ denotes the exceptional divisor 
in $\widetilde{{\Bbb C}^2}$ and that the 
Hirzebruch surface $\Sigma_1 = {\Bbb P}({\cal O}(-1) \oplus {\cal O})$
is the complex compactification of 
$\widetilde{{\Bbb C}^2}$ obtained by adding a line at infinity 
$\ell_{\infty}.$ The following lemma is easy to prove.

\begin{lemma}\label{sigma}There is a bijection between the sets
${\cal E}_2(\widetilde{{\Bbb C}^2},\widetilde{{\Bbb C}^2}- \ell)$
and ${\cal E}_2(\Sigma_1, \ell_{\infty}).$  
\end{lemma}
 
The proof is in the appendix.

\subsection{Moduli spaces}
We wish to study moduli  of instantons. 
We  say that  two triples $\Delta= (j,p,t_{\infty})$
 and $\Delta'=(j',p',t'_{\infty})$ 
are equivalent if they represent the
same instanton. 
In terms of holomorphic bundles, this
means that two triples are equivalent 
if their corresponding holomorphic bundles $E$ and $E'$ over 
 $\widetilde{{\Bbb C}^2}$ framed at infinity  
are isomorphic, via an isomorphism taking 
$t_{\infty}$ into $t'_{\infty}.$
In particular 
these bundles give isomorphic restrictions over the exceptional divisor, 
hence $E$ and $E'$
 must have the same splitting type, that is, $j=j'.$ 

Let us  consider triples 
$(j,p,t_{\infty}) $ and $(j,p',t'_{\infty}), $
with the same integer $j,$ 
 representing holomorphic bundles 
$E$ and $E'$ over  $\widetilde{{\Bbb C}^2}$
trivialized at infinity.
We know
from proposition \ref{sigma} that these bundles 
 may be looked upon as 
bundles $(E,t_{\infty})$ 
and $(E',t'_{\infty})$
over $\Sigma_1$
trivialized over ${\ell}_{\infty}.$
   An isomorphism for  framed bundles 
is a bundle isomorphism $\Phi \colon E \rightarrow E'$
such that $\Phi(t_{\infty})=t'_{\infty}.$
Two framings  $t_{\infty}$ and $ t'_{\infty}$ for the same   
same underling bundle $E$ over $\Sigma_1$
differ by a holomorphic map 
$\Phi\colon  {\ell}_{\infty}\rightarrow {\mbox GL}(2,{\Bbb C})$
and,  since $\ell_{\infty}$ is compact, $\Phi$ must be constant.
Hence, projecting  $(E,t_{\infty})$ 
on the first coordinate we obtain  a fibration of 
the space of framed bundles over $\Sigma_1$ 
 over the space of bundles over $\Sigma_1$ which are trivial 
on the line at infinity, with fibre ${\mbox GL}(2,{\Bbb C}).$

\begin{equation}\label{fibration} \begin{array}{c}
{\mbox GL}(2,{\Bbb C}) \\
\downarrow \\
\left\{ framed \, rank-2 \, bundles \, over \, \Sigma_1 \right\} \\
   \downarrow \\
\left\{ rank-2 \, bundles\, over \, \Sigma_1 \, trivial \,on \, {\ell}_{\infty} \right\}.
\end{array}\end{equation}
\vspace{3mm}

We
 are thus led to  study the base space of
this fibration, or equivalently, the space of isomorphism classes 
of bundles on $\widetilde{{\Bbb C}^2}$
which are trivial on the line at infinity.
We  define 
${\cal M}_j$ to be   space
of rank two holomorphic  bundles on the $\widetilde{{\Bbb C}^2}$
with vanishing first Chern class and with splitting type $j,$ modulo
isomorphism, that is,
\begin{equation}\label{Mj}
{\cal M}_j = \left.\left\{ \begin{array}{ll} E \,\ hol. \,\, bundle \,\,  
 over \,\, \widetilde{{\Bbb C}^2}: \\
E|_{\ell} \simeq  
{\cal O} (j) \oplus {\cal O}(-j) \end{array} \right\}\right/ \sim. \end{equation}
\vspace{2mm}

Fix the splitting type $j$ and set $J=(j-1)(2j-1),$ then the polynomial
$p$ has $J$ coefficients. 
Identifying  the polynomial $p$ with the $J-$tuple 
formed by its coefficients written in lexicographical order, we
may define 
 in ${\Bbb C}^J$ 
the equivalence relation $p \sim p'$ 
if $(j,p)$ and $(j, p')$ represent isomorphic
bundles.
Set-theoretically there is an identification
\begin{equation}\label{quotient}  {\cal M}_j
   = {\Bbb C}^J/\sim.\end{equation}
We give ${\Bbb C}^J/\sim$
the quotient topology and  ${\cal M}_j$ 
the topology induced by (\ref{quotient}). 
${\cal M}_j$ 
is generically a complex projective space     
of dimension $2j-3$ ( \cite{JA} Thm.\thinspace3.5).
However the topology of ${\cal M}_j$ is quite complex, and, 
in particular,  is non-Hausdorff for any $j\geq 2.$

There is a topological   embedding taking  
${\cal M}_j$  into the least generic strata of ${\cal M}_{j+1}.$
We write it out explicitly in 
coordinates, representing an element  $E \in {\cal M}_j,$
by its  canonical form of transition matrix 
according to (\ref{matrix}).

\begin{guess}: The following map defines a topological embedding 
$$\begin{array}{rcl}\Phi_j: {\cal M}_j & \rightarrow &  {\cal M}_{j+1} \cr
 (j,p) & \mapsto &  (j+1,zu^2p) \end{array}. $$
\end{guess}

The proof is in the appendix.

The map $\Phi_j$ takes ${\cal M}_j$ into the least generic strata 
of ${\cal M}_{j+1}.$ In fact, im$\Phi$ is the subset of 
${\cal M}_{j+1}$ consisting of bundles that split in the second
 formal
neighborhood of the exceptional divisor.
The complexity of the  topology of ${\cal M}_j$ increases with $j$
according to these embeddings. ${\cal M}_2$ is non-Hausdorff 
and therefore this property persists in ${\cal M}_j$ for $j \geq2.$
Explicitly ${\cal M}_2 \simeq {\Bbb P}^1 \cup \{A,B\}$ where
the generic set ${\Bbb C}P^1$ consists of bundles that do not 
split on the first formal neighborhood, $A$ and $B$ 
are special points corresponding to two special  bundles with splitting
type $2;$ the one  
that splits on the first
formal neighborhood but not on higher neighborhoods, and the split bundle 
(see \cite{JA})

\subsection{Instanton numbers}

We now define the instanton numbers that stratify the 
spaces ${\cal M}_j$ into Hausdorff components. 

We consider a compact  complex (smooth) surface  $X$ 
together with the blow--up 
 $\pi
\colon \widetilde{X} \rightarrow
X$ of a point $x \in X$  and  once again denote by 
$\ell$ the 
exceptional divisor.
Let $\widetilde{E}$
be a rank 2 holomorphic bundle 
  over 
$\widetilde{X}$ satisfying
$\det \widetilde{E} \simeq {\cal O}_ {\widetilde{X}}.$ 
The splitting type of $\widetilde{E}$ 
is by definition the
integer $j \geq 0$ such that 
$\widetilde{E}|_{\ell}  \simeq {\cal O}(j) \oplus
 {\cal O}(-j).$
Set  $E = \pi_* {\widetilde{E}}^{\vee \vee}.$
Assuming  $X$ compact,
Friedman and Morgan  \cite{FM}, p.\thinspace393
 gave the following estimate 
relating the second Chern classes
to the splitting type
$$  j \leq  c_2(\widetilde{E}) - c_2(E) \leq j^2.$$
Sharpness of these bounds was proven in   \cite{CA2}.
Let $F$ be
a bundle on $\widetilde{{\Bbb C}^2}$
with vanishing first Chern class.
If $X$ is a compact complex surface, then there exist holomorphic bundles
$\widetilde{E} \rightarrow \widetilde{X}$ 
which are isomorphic to $F$ on a 
neighborhood of the exceptional divisor.
In fact, following \cite{BSPM}, given a bundle $E \rightarrow X,$
we can construct 
bundles 
 $\widetilde{E} \rightarrow \widetilde{X}$
 satisfying:

\noindent$\iota )$ 
$\widetilde{E}|_{\widetilde{X}-{\ell}} = \pi^*(E|_{X-p}),$
where $\pi$ is the blow--up map, and

\noindent$\iota\iota)$ $\widetilde{E}|_{V} \simeq F|_{U}$  for small 
 neighborhoods $V$ and $U$ of the exceptional divisor in 
$\widetilde{X}$ and $\widetilde{{\Bbb C}^2}$ respectively.

Moreover, every bundle $\widetilde{E}$
on $\widetilde{X}$ is obtained 
 this way  [8, Cor.\thinspace3.4].
The isomorphism class of $\widetilde{E}$ depends on the attaching map
$\phi\colon (\widetilde{X}-\ell) \cap V \rightarrow {\mbox GL}(2,{\Bbb C }),$
however, the topological type of $\widetilde{E}$ is independent of 
$\phi.$ 
 Therefore the charge
does not depend upon the choice of $\phi.$
Since the $n-$th infinitesimal neighborhood 
of the exceptional divisor on a compact complex surface
is isomorphic (as a scheme) to the $n-$th infinitesimal 
neighborhood of the exceptional divisor on $\widetilde{{\Bbb C}^2}$
we are able to use an explicit description 
 for bundles on $\widetilde{{\Bbb C}^2},$  even though 
$x \in \widetilde{X}$ might not have an 
open neighborhood analytically equivalent 
to $\widetilde{{\Bbb C}^2}.$ We quote:

\begin{guess}\label{gluing} (\cite{BSPM}, Cor.\thinspace4.1)
Let $X$ be a compact surface and $\widetilde{X}$ 
denote the blow up of $X$ at $x.$
Every
holomorphic rank 2 vector $\widetilde{E}$ bundle over $\widetilde{X}$
with vanishing first Chern class
is topologically  determined
 by a triple $(E,j, p)$
where $E$ is a rank 2 holomorphic  bundle on $X$ with
vanishing first Chern class, $j$ is a nonnegative
integer, and  $p$ is a polynomial.
\end{guess}

If $\widetilde{E}$ is as in the above proposition, we denote 
\begin{equation} \label{(j,p)} \widetilde{E}:=(E,j, p). \end{equation}
The pair $(j, p)$ gives an explicit description of $\widetilde{E}$
on a neighborhood of the exceptional divisor, and determines 
the charge of $\widetilde{E}.$ 
To
calculate
the charge, we actually compute
 two finer numerical invariants of $\widetilde{E},$
 which we now describe.
Following Friedman and Morgan (\cite{FM}, p.\thinspace302), we
 define a sheaf  $Q$ by the
exact sequence,
$$
0 \rightarrow \pi_* \widetilde{E}
 \rightarrow \pi_*(\widetilde{E})^{\vee\vee} \rightarrow Q \rightarrow 0.
$$
Note that $Q$ is supported only at the point $x.$
 From the exact sequence $(1.6)$
it
follows immediately that
 $ c_2(\pi_*\widetilde{E}) - c_2(E) = l(Q),$
where $l$ stands for length.
An application of
Grothendieck--Riemann--Roch (see \cite{HA}, p.\thinspace392) gives
that $$c_2( \widetilde{E}) - c_2(E) = l(Q) +
l(R^1 \pi_* \widetilde{E}).$$
We call $w\colon =l(Q) $ the {\em width} 
and 
$h \colon = l(R^1 \pi_* \widetilde{E})$
the ${\em height}$ 
of the instanton $\widetilde{E}.$ 

\subsection{Holomorphic instanton patching}
Let $X$ be a surface with polarization $L.$
Choose $N>>0,$ then
$\widetilde{L} = NL - {\ell}$ is ample 
and  it is natural to choose $\widetilde{L}$ as a polarization
of $\widetilde{X}.$ We fix these choices of polarizations and 
by stable bundle we mean stable with respect to the fixed polarization.
By the Hitchin--Kobayashi correspondence for compact surfaces, 
instantons 
 correspond to stable bundles. 
A stable bundle $\widetilde{E}$ on $\widetilde{X}$ 
such that 
$\widetilde{E}\vert_{\widetilde{X} -{\ell}} \simeq \pi^*(E\vert_{X -\{p\}})$
with $E$ stable on $X$ 
and $\widetilde{E}\vert_{V({\ell})} \simeq  \Delta=(j,p,t_{\infty})$ 
on some neighborhood $V({\ell})$ of the exceptional divisor
is said to be obtained by obtained by {\em holomorphic patching} 
of the instantons of $\Delta$ to $E.$ The reason for this terminology is 
that given $E$ and $\Delta$ any choice of
gluing $\Phi : {\Bbb C}^2 -\{0\} \rightarrow {\mbox GL}(2,{\Bbb C})$
gives  a holomorphic way to construct a new 
instanton.   Equivalently, it is enough to 
choose a framing at the point $p.$

\begin{lemma} Every instanton on $\widetilde{X}$ is 
obtained by {\em ``holomorphic''} patching of an 
instanton on $\widetilde{{\Bbb C}^2}$  to an instanton on $X.$
\end{lemma}
 
\noindent{\bf Proof}:
By \cite{BSPM} Corollary\thinspace 3.4 
every holomorphic rank two vector bundle 
$\widetilde{E}$ over $\widetilde{X}$
with vanishing first Chern class 
is completely determined (up 
to isomorphism) by a 4-tuple $\widetilde{E}\colon = (E,j,p,\Phi)$ 
where $E$ is a rank two holomorphic  bundle on $X$ with 
vanishing first Chern class, $j$ is a nonnegative 
integer, $p$ is a polynomial, 
and $\Phi : {\Bbb C}^2 -\{0\} \rightarrow {\mbox GL}(2,{\Bbb C})$ is
a holomorphic map. The bundle $\widetilde{E}$
has splitting type $j$ over the exceptional divisor,
 and satisfy the property
$\widetilde{E}\vert_{\widetilde{X} -{\ell}} \simeq \pi^*(E\vert_{X -\{p\}}).$
If $\widetilde{E}$ is stable, then so is $E$ (see \cite{FM}).\hfill\square{5}  

\begin{remark} The charge addition given by the  patching of 
$\Delta$ can be calculated by a Macaulay2 program written by 
Irena Swanson and the author \cite{M2}. The program has as input 
 $j$ and $p$ and as output the height and the width of an
instanton $(j,p).$ 
\end{remark}

\subsection{Stratification of Moduli of Instantons}

The following theorem shows that instanton numbers provide good stratifications
for moduli of instantons on $\widetilde{\Bbb C}^2.$ In fact, these 
numbers stratify the spaces ${\cal M}_j$ into Hausdorff components, 
and this is the coarsest stratification of 
${\cal M}_j$ for which the strata are Hausdorff. In \cite{PAMS} 
it is shown that stratification by Chern numbers is not
fine enough to have this property. We cite.

\begin{theorem}(\cite{PAMS} Thm.\thinspace 4.1)
\label{stratification} The numerical invariants $w$
and $h$ provide a  decomposition ${\cal M}_j
= \cup S_i$ 
where each $S_i$ is homeomorphic to an open subset
of a complex projective space of dimension at most $2j-3.$
 The lower bounds for these invariants are $(1,j-1)$ and this pair of
invariants takes
  place on the generic part of ${\cal M}_j$
which is homeomorphic to   ${\Bbb CP}^{2j-3}$
minus a closed subvariety of
codimension at least 2.
The   upper bounds for these invariants are $(j(j-1)/2, j(j+1)/2)$
and this pair occurs at one single point of ${\cal M}_j$
which represents the split bundle.
\end{theorem}

\section{Curve singularities}

 Here is how to use
instanton numbers 
  to distinguish curve singularities.
Start with a curve $p\,(x,y)=0$ on
${\Bbb C}^2.$ Choose your favorite integer $j$ 
and construct an instanton on   $\widetilde{{\Bbb C}^2}$i
having data $(j,p).$ Calculate the  height and the width of 
the instanton, use them as analytic invariants
of the curve, and use the charge as a topological invariant.
In other words, we are using the polynomial 
defining the plane curve as an extension class in 
$\mbox{Ext}^1{\cal O}(j), \cal O(-j)). $ This defines a bundle 
$E(j,p)$ as in \ref{(j,p)}. We then calculate  the instanton 
numbers of this bundle, as defined in section 2,  and 
regard them as being associated to the curve.

Note that 
to perform the computations  we must choose a  representative 
for the curve and coordinates for the bundle. 
I
use the canonical choice of coordinates for $\widetilde{{\Bbb C}^2}$
 as in section 2.
Taking into account that the blow-up map in these 
coordinates is given by $x \mapsto u$ and $y \mapsto zu$ 
the bundle $E(j,p)$ is then given canonically in these coordinates by
$$E(j,p)\colon =\left(\matrix{ z^j & p(u,zu) \cr 0 &  z^{-j} \cr} \right).
$$
If a second representative 
$\bar p$ for the same singularity is given,  there is a holomorphic 
change of coordinates $\phi$ taking $p$ to $\bar p.$
To compute the invariants using this second representative, the 
coordinate change has to be applied to the bundle as well. 
In this paper I give only a couple of results to illustrate 
the behavior of the instanton numbers applied to singularities. 
Explicit hand-made computations of these invariants for 
small values of $j$ appear in \cite{PAMS} and \cite{CA2}.
The invariants can be computed by a Macaulay2 algorithm
written by Irena Swanson and the author, see Remark \ref{program}.

The next theorems
show that instanton numbers distinguish the most basic singularities
and also give some examples where instanton numbers are finer than 
classical invariants.

\begin{theorem} Instanton numbers distinguish nodes/tacnodes
from  cusps/higher order cusps.
\end{theorem}

\noindent{\bf Proof}: These singularities have  quasi--homogeneous 
representatives of the form $y^n-x^m,$ $n < m,$ 
$n$ even for nodes and tacnodes, and
 $n$ odd for cusps and higher order cusps.
 We want to show that instanton numbers detect 
 the parity of the smallest these exponents. In fact, more is true, 
instanton numbers detect the multiplicity itself.

Suppose $n_1 < n_2.$ We claim that if  $j >n_2$ then  $w(j,p_1) \neq 
w(j,p_2).$ In fact, for $n<m$ and large enough $j$ 
the width takes the value  $$w(j,y^n-x^m)=   n(n+1)/2.$$ 
Alternatively, by vector bundle reasons we have that
 $w(j,p_1) < w(j,p_2).$ The second assertion is easier 
to show. The holomorphic bundle 
$E(j,p_1)$ restricts as a non-trivial extension on  
the $n_1$th formal  neighborhood  $l_{n_1 }$   
whereas  $E(j,p_2)$ splits on  $l_{n_1 }.$ 
These bundle therefore  belong to different strata of ${\cal M}_j$
and by theorem \ref{stratification} must have
distinct instanton numbers.  \hfill\square{5}

\begin{theorem} In some cases instanton numbers are finer than classical
invariants.
\end{theorem}

\noindent{\bf Proof}: See   tables I and II below. \hfill\square{5}

\vspace{5mm}
The classical invariants we consider are:
\begin{itemize}
\item $\delta_P = dim(\widetilde{\cal O}_P/{\cal O}_P) $ 
\item Milnor number $\mu = {\cal O}/ <J(P)>$ 
\item Tjurina number $\tau = {\cal O}/<P,J(P)>$
\end{itemize}

\vspace{5mm }
\noindent{\bf Note}: The 
first table is motivated by
 exercise 3.8 of Hartshorne \cite{HA} page 395. However,
in the statement of the problem, 
the first polynomial contains an incorrect exponent. It is  written as 
``$x^4y-y^4$''
but it should be ``$x^5y-y^4.$''
 
\vspace{5mm}
\begin{tabular}{|l|l|l|l|l|l|}
\hline
 TABLE I & & & &  \multicolumn{2}{|c|}{  $j=4$}  \\     
\hline
polynomial         &     $\delta_P$ &  $\mu$ &  $\tau$ &   $w$ &  $h$\\
\hline
\hline
$ x^5y-y^4$             &     9 &          17 &    17 &     10 &   6\\
\hline
$ x^8-x^5y^2-x^3y^2+y^4$  &    9  &         17  &   15  &     8  &   6\\
\hline
\end{tabular}

\vspace{5mm}
\noindent The second table shows an example where instanton numbers are
finer than $\delta_P, \mu, $ and $\tau.$

\vspace{5mm}
\begin{tabular}{|l|l|l|l|l|l|}
\hline
 TABLE II & & & &  \multicolumn{2}{|c|}{  $j=4$}  \\     
\hline
polynomial         &   $\delta_P$ & $\mu$ &  $\tau$&   $w$&  $h$\\
\hline
\hline
 $x^2-y^7$           &     3 &             6  &  6 &      3 &    5\\
\hline
 $x^3-y^4$           &     3    &          6  &     6 &    6   &   6\\
\hline
\end{tabular}

\vspace{5mm}

\begin{remark} The idea of using the polynomial defining a
 singularity as the
 extension class of a holomorphic bundle can be further
 generalized in several ways. For curves themselves, one can use 
other base spaces. For instance, constructing bundles on the 
total space of ${\cal  O}_{{\Bbb P}^1}(-k)$ requires very little modifications,
but give quite different results. 
One can also generalize to hypersurfaces in higher dimensions. 
\end{remark}

\section{Appendix}

This appendix contains  proofs of two technical but straightforward
results used in the text.

\vspace{5mm}
\noindent{\bf Proof of Lemma 2.1}:
 Given that $\Sigma_1 = \widetilde{{\Bbb C}^2}
\cup \ell_{\infty},$ there exists an open  neighborhood $W$ 
of $\ell_{\infty}$ satisfying:

$\iota)$ $\Sigma_1 - W = \ell$

$\iota\iota)$ $W - \ell_{\infty} \simeq  \widetilde{{\Bbb C}^2} - \ell$

$\iota\iota\iota)$ we have a commutative diagram 

$$\begin{matrix}{ \ell_{\infty} & \rightarrow & {\Bbb P}^1 \cr
                          \downarrow & & \downarrow \cr
                    W & \rightarrow & {\cal O}(1)}\end{matrix}$$
where the vertical arrows are inclusions and the horizontal 
arrows are isomorphisms.
Now, given $(E,\eta) \in 
{\cal E}_2(\widetilde{{\Bbb C}^2},\widetilde{{\Bbb C}^2}- \ell),$
i.e., $\eta = (a,b);$ where $a,b\colon \widetilde{{\Bbb C}^2}- \ell
\rightarrow E|_{\widetilde{{\Bbb C}^2}- \ell}$ are 
linearly independent sections, define
$$\Phi(E,\eta):= (F,\mu) \in {\cal E}_2(\Sigma_1,\ell_{\infty})$$
by gluing $E$ with 
$W \times {\Bbb C}^2$ along $\eta,$ that is, define
$$F= E \sqcup (W\times{\Bbb C}^2)/ \sim$$
where
$\alpha\, a(x) + \beta\, b(x) \sim (x,(\alpha,\beta))$
for  $x \in \widetilde{{\Bbb C}^2}- \ell$ and $(\alpha,\beta)
\in {\Bbb C}^2.$

Let $\mu= (e_1,e_2)$ be the canonical section of 
the trivial bundle over $\ell_{\infty},$ that is,
$\mu(y)= (y,(1,0),(0,1))$ where $y \in \ell_{\infty} \subset W.$ 
Now, given $(F,\mu) \in {\cal E}_2(\Sigma_1,\ell_{\infty})$
define
$$\Psi(F,\mu):=(E,\eta)$$
where $E=F|_{\widetilde{{\Bbb C}^2}}$
and $\eta$ is defined as follows.
Since $F|_W$ is trivial (because 
$F\left|_{W|_{\ell_{\infty}}}\right. = F|_{{\ell}_{\infty}}$
is trivial) there is a unique trivialization 
$\widetilde{\eta}$ over $W$ such that 
$\widetilde{\eta}|_{{\ell}_{\infty}}= \mu$
(this is because $H^0(W,{\cal O})=H^0({\cal O}(1),{\cal O})={\Bbb C}).$
Define $$\eta:=\widetilde{\eta}|_{W_{{\ell}_{\infty}}}= 
\widetilde{\eta}|_{\widetilde{\Bbb C}-\{0\}}.$$

It is straightforward to prove that 
$\Psi\Phi$ and $\Phi\Psi$ are the identities.\hfill\square{5}

\vspace{5mm}
\noindent{\bf Proof of Proposition 2.2:}
We first show that the map is well defined. 
Suppose 
$\left(\matrix{ z^j & p \cr 0 & z^{-j} }\right)$ and 
$\left(\matrix{ z^j & p' \cr 0 & z^{-j} }\right)$ 
represent isomorphic bundles. 
Then there are coordinate changes
 $ \left(\matrix{a  & b   \cr c  & d \cr}\right)$
holomorphic in $z, \, u$ and 
$ \left(\matrix{\alpha & \beta  \cr \gamma & \delta \cr}\right)$
  holomorphic in $z^{-1}, \, zu$
for which the following equality holds
(compare \cite{CA1} pg. 587)
$$ \left(\matrix{\alpha & \beta  \cr \gamma & \delta \cr}\right)
 = \left(\matrix{z^j & p'  \cr 0 & z^{-j} \cr}\right)
 \left(\matrix{a  & b   \cr c  & d \cr}\right)  
\left(\matrix{z^{-j} & -p  \cr 0 & z^j \cr}\right).$$
Therefore these two bundles are isomorphic 
exactly when the system of equations 
$$ \left(\matrix{\alpha & \beta  \cr \gamma & \delta \cr}\right)
= \left(\matrix{a + z^{-j}p'c & 
z^{2j}b +z^j(p'd  -ap) - pp'c  \cr
z^{-2j}c & d - z^{-j}pc \cr}\right) \eqno (*)$$
can be solved by a matrix 
 $ \left(\matrix{a  & b   \cr c  & d \cr}\right)$
 holomorphic in $z, \, u$ which makes 
$ \left(\matrix{\alpha & \beta  \cr \gamma & \delta \cr}\right)$
 holomorphic in $z^{-1}, \, zu.$

On the other hand, the images of these two bundles are given by
transition matrices
$\left(\matrix{ z^{j+1} & z\,u^2p \cr 0 & z^{-j-1} }\right)$ and 
$\left(\matrix{ z^{j+1} & z\,u^2p' \cr 0 & z^{-j-1} }\right),$
which  represent isomorphic bundles iff there are
coordinate changes
  $ \left(\matrix{\bar a  & \bar b   \cr \bar c  & \bar d \cr}\right)$
holomorphic in $z, \, u$ and 
$ \left(\matrix{\bar\alpha & \bar\beta  \cr \bar\gamma & 
\bar\delta \cr}\right)$
  holomorphic in $z^{-1}, \, zu$ 
satisfying the equality
$$ \left(\matrix{\bar\alpha & \bar\beta  \cr \bar\gamma & 
\bar\delta \cr}\right)
 = \left(\matrix{z^{j+1} & z\,u^2p'  \cr 0 & z^{-j-1} \cr}\right)
 \left(\matrix{\bar a  & \bar b   \cr \bar c  & \bar d \cr}\right)  
\left(\matrix{z^{-j-1} & -z\,u^2p  \cr 0 & z^{j+1} \cr}\right).$$
That is, the images represent isomorphic bundles if the system
$$ \left(\matrix{\bar\alpha & \bar\beta  \cr \bar\gamma & 
\bar\delta \cr}\right)
= \left(\matrix{\bar a + z^{-j } u^2p'\bar c & 
z^{2j+2}\bar b +z^{j+2}u^2(p'\bar d  -\bar a p) - z^2u^4pp'\bar c  \cr
z^{-2j-2}\bar c & \bar d - z^{-j} u^2p\bar c \cr}\right) \eqno (**)$$
has a solution.

 Write $x = \sum x_i u^i$ for $x \in \{a,b,c,d,
 \bar a,\bar b,\bar c,\bar d\}$
and choose
$\bar a_i = a_{i+2},$ 
$\bar b_i = b_{i+2}u^2,$
$\bar c_i = c_{i+2}u^{-2},$
$\bar d_i = d_{i+2}.$  
Then if $ \left(\matrix{a  & b   \cr c  & d \cr}\right)$ 
solves  (*), one verifies  that 
 $ \left(\matrix{ \bar a  & \bar b   \cr \bar c  & \bar d \cr }\right)$
solves 
(**), which implies that the images  represent 
 isomorphic bundles and therefore $ \Phi_j $  is well defined.
To show that the map is injective just reverse the 
previous argument.
Continuity is obvious. 
Now we observe also that the image 
 $\Phi_j({\cal M}_j)$ is a saturated set in ${\cal M}_{j+1}$
 (meaning that if 
$ y \sim x$ and $x \in \Phi_j({\cal M}_j)$ then $y \in \Phi_j({\cal M}_j$)).
In fact, if $E \in 
 \Phi_j({\cal M}_j)$ 
then $E$ 
splits in the 2nd formal  neighborhood.
Now if $E' \sim E$ than $E'$ must also
split in the 2nd formal  neighborhood
therefore the polynomial corresponding to $E'$ is 
of the form $u^2p'$
and hence $ \Phi_j(z^{-1}p') $ gives $E'.$
Note also that 
 $\Phi_j({\cal M}_j)$ is a closed 
subset of ${\cal M}_{j+1},$ given by
the equations 
$p_{il} = 0$ for $i = 1,2$ and  $i-j+1 \leq l \leq j-1.$
Now the fact that  $\Phi_j$ is a homeomorphism over its  
image follows from the following  easy lemma.\hfill\square{5}

\begin{lemma} 
Let $X \subset Y$ be a closed subset 
and $\sim$ an equivalence relation in $Y,$ such 
that $X$ is $\sim$ saturated. 
Then the map $I :X/{\sim} \rightarrow Y/{\sim}$
induced by the inclusion is a homeomorphism over the image.
\end{lemma}

\noindent{\bf Proof:} Denote by 
$\pi_X : X \rightarrow X/{\sim} $
and $\pi_Y: Y \rightarrow Y/{\sim}$ the projections.
Let $F$ be a closed subset of $X/{\sim}.$ Then 
$\pi_X^{-1}(F)$ is closed and saturated in $X$ and therefore 
$\pi_X^{-1} (F) $ is also closed and saturated in $Y.$ 
It follows that $\pi_Y(\pi_X^{-1}(F))$  is closed 
in $Y/{\sim}.$ \hfill\square{5}

\vspace{5mm}

\noindent{\small Elizabeth Gasparim}

\noindent{\small New Mexico State University} 

\noindent{\small Department of Mathematics}

\noindent{\small Las Cruces NM 88001}

\noindent{\small gasparim@nmsu.edu}

\end{document}